\documentclass[12pt]{article}
\usepackage{amscd, amssymb, latexsym, amsmath}

\newtheorem{thm}{Theorem}

\newtheorem{lem}{Lemma}

\newenvironment{proof}
{\noindent {\em Proof.}} {\hfill $\Box$}

\numberwithin{thm}{section} \numberwithin{cor}{section}
\numberwithin{pro}{section}
\numberwithin{lem}{section} \numberwithin{equation}{section}

\newcommand{\R}{\mathbb R}

\numberwithin{equation}{section}

\newcommand{\pai}{\frac{\partial}{\partial x^i}}

\begin{document}
\title
{Stability and curvature estimates for minimal graphs with flat
normal bundles}
\author{Mu-Tao Wang}
\date{September, 2004, revised November 7, 2004}
\maketitle
\begin{abstract}
It is well-known that a minimal graph of codimension one is
stable, i.e. the second variation of the area functional is
non-negative. This is no longer true for higher codimensional
minimal graphs in view of an example of Lawson and Osserman. In
this note, we prove that a minimal graph of any codimension is
stable if its normal bundle is flat.  We also prove minimal graphs
of dimension no greater than six and any codimension is flat if
the the normal bundle is flat and the density at infinity is
finite.  Such a Bernstein type theorem holds in any dimension if
we assume additionally growth conditions on the volume element.

\end{abstract}

\section{Introduction}

The graph of a solution $f:D\subset \R^n\rightarrow \R$ of the
minimal surface equation

\begin{equation}\label{mse}
div(\frac{\nabla f}{\sqrt{1+|\nabla f|^2}})=0
\end{equation} is naturally a minimal hypersurface in $\R^{n+1}$.
In general, we may consider the graph of a vector-valued function
and ask when this is a minimal submanifold of the Euclidean space.
The function then satisfies a nonlinear elliptic system. Indeed, a
$C^2$ vector-valued
 function $f=(f^1, \cdots f^m):D
\rightarrow \R^m$ is said to be a solution to the minimal surface
system (see Osserman \cite {os1} or Lawson-Osserman \cite{lo}) if

\begin{equation}\label{mss}
\sum_{i,j=1}^n\pai(\sqrt{g}g^{ij} \frac{\partial
f^\alpha}{\partial x^j})=0\,\,\text{for each } \alpha=1\cdots m
\end{equation} where $ g^{ij}=(g_{ij})^{-1}$,
$g_{ij}=\delta_{ij}+\sum_{\beta=1}^m \frac{\partial
f^\beta}{\partial x^i}\frac{\partial f^\beta}{\partial x^j}$, and
$g=\det g_{ij}$. The graph of $f$, so called a {\it minimal
graph}, is then a minimal submanifold of $\R^{n+m}$ of dimension
$n$ and codimension $m$.

 The minimal surface system was first studied in
Osserman \cite{os1}, \cite{os2} and Lawson-Osserman \cite{lo}. In
contrast to the codimension one case, Lawson and Osserman
\cite{lo} discovered remarkable counterexamples to the existence,
uniqueness and regularity of solutions to the minimal surface
system in higher codimension. It is thus interesting to identify
natural conditions under which  theorems for the minimal surface
equation can be generalized.

The difficulty of the higher codimensional problems is amplified
by the complexity of the normal bundle. Given an $n$ dimensional
submanifold $\Sigma$ of $\R^{n+m}$, recall the normal bundle
consists of the orthogonal complement of the tangent spaces of
$\Sigma$ in $\R^{n+m}$. Near a point of $\Sigma$, choose an
orthonormal frame $e_1, \cdots, e_n$ for the tangent bundle and
$e_{n+1},\cdots e_{n+m}$ for the normal bundle. The coefficients
of the second fundamental form is denoted by $h_{\alpha
ij}=\langle \nabla_{e_i} e_\alpha, e_j\rangle$. Recall from the
Ricci equation, the curvature of the normal bundle is given by

\begin{equation}\label{normal}R_{\alpha \beta ij}=h_{\alpha ik}h_{\beta jk}-h_{\alpha
jk}h_{\beta ik}.\end{equation}

We say $\Sigma$ has {\it flat normal bundle} if $R_{\alpha \beta
ij}\equiv 0$, see for example \cite{te}. When $\Sigma$ is a graph,
we can in fact choose $e_{n+1}, \cdots e_{n+m}$ to be globally
parallel sections. Equation (\ref{normal}) holds trivially when
$\Sigma$ is of codimension one, i.e. $m=1$. In particular, any
oriented hypersurface has flat normal bundle.

Recall a minimal submanifold is called {\it stable} if the second
variation of the volume functional is non-negative with respect to
any compact-supported variation fields. By a calibration argument,
a minimal graph of codimension one is stable (in fact
area-minimizing). This is no longer true in higher codimension by
a counterexample of Lawson and Osserman \cite{lo}. Nevertheless,
we prove the following stability theorem.

\begin{thm} If $\Sigma $ is a minimal graph with flat normal bundle
in the Euclidean space, then $\Sigma$ is stable.\end{thm}

 A different stability criterion for higher codimensional
minimal graphs is studied in \cite{lw}. We first generalize the
curvature estimate of Schoen-Simon-Yau \cite{ssy} and prove the
following Bernstein type theorem.

\begin{thm}
Suppose $\Sigma $ is the graph of an entire smooth function
$f:\R^n\rightarrow \R^m$ of the minimal surface system for $n\leq
5$. If the normal bundle of $\Sigma$ is flat and $Vol(\Sigma \cap
B_R)\leq c R^{n}$ for some constant $c$, then $f$ is a linear map.
\end{thm}
Here $B_R$ is the ball of radius $R$ in $\R^{n+m}$ centered at the
origin.

Ecker and Huisken \cite{eh} prove the following Bernstein type
result in the codimension one case.

\noindent {\bf Theorem}\hskip 5pt {\it Suppose $\Sigma$ is the
graph of an entire smooth solution $f:\R^n\rightarrow \R$ of the
minimal surface equation. $f$ is a linear map if
\[\sqrt{1+|Df|^2}(x)=o(\sqrt{|x|^2+|f(x)|^2}).\]}

 We  generalize Ecker and Huisken's theorem \cite{eh} to the higher codimensional case.

\begin{thm}
Suppose $\Sigma$ is the graph of  an entire smooth solution
$f:\R^n\rightarrow \R^m$of the minimal surface system. $f$ is a
linear map if the following three conditions hold:

\indent 1) the normal bundle of $\Sigma$ is flat.

\indent 2) \[\sqrt{\det(I+(df)^T
df)}(x)=o(\sqrt{|x|^2+|f(x)|^2}).\]

\indent 3) $Vol(\Sigma\cap B_R)\leq c(n) R^n$ for some constant
$c(n)$.
\end{thm}

Higher codimensional Bernstein type theorems have been studied by
many authors \cite{fc}, \cite{hjw}, \cite{jx}, \cite{wa3} assuming
various conditions.

The author would like to thank Professor C.-L. Terng for
suggesting him to look at mean curvature flows of submanifolds of
flat normal bundles in the spring of 2002. This note was written
up while the author was teaching a graduate course at Columbia in
which he went over curvature estimates for minimal hypersurfaces.
During when he realized the key formulae (see equations
(\ref{eq1}) and (\ref{eq2})) to generalize Schoen-Simon-Yau
\cite{ssy} and Ecker-Huisken\cite{eh} to the flat normal bundle
case were contained in his earlier work \cite{wa1} and \cite{wa2}.
With these formulae, the derivations follow straightforward from
\cite{ssy} and \cite{eh}.

We remark that mean curvature flows of submanifolds with flat
normal bundles are studied in a recent paper by Smoczyk, Wang and
Xin \cite{swx}.

\section{Preliminary}

Let $\Sigma$ be the graph of $f:D\subset \R^n\rightarrow \R^m$ and
$\{e_1,\cdots, e_n\}$ and $\{e_{n+1}, \cdots e_{n+m}\}$ be local
orthonormal bases for $T\Sigma$ and $N\Sigma$, respectively.

Let $\Omega$ denote the n-form $dx^1\wedge\cdots\wedge dx^n$ where
$x^1\cdots, x^n$ are coordinates on $\R^n$. We extend $\Omega$ to
$\R^{n+m}$ and consider the function  $*\Omega=\Omega(e_1, \cdots,
e_n)$ on $\Sigma$. Notice $*\Omega$ is the Jacobian of the natural
projection from $\Sigma$ to $D$ and $*\Omega>0$ on $\Sigma$. In
terms of $f$, we have  $*\Omega=\frac{1}{\sqrt{\det(I+(df)^T
df)}}$.

 We recall the following formula derived in
\cite{wa2}(Proposition 3.1) and \cite{wa3}
\begin{equation}
\Delta*\Omega +*\Omega |A|^2+2\sum_{k=1}^n\sum_{\alpha, \beta, i<j
}\Omega_{\alpha\beta ij} h_{\alpha i k}h_{\beta jk}=0\\
\end{equation}
where $\Omega_{\alpha\beta ij}=\Omega (e_1,\cdots,
e_\alpha,\cdots, e_\beta,\cdots, e_n)$ with $e_\alpha$ occupying
the $i$-th place and $e_\beta$ occupying the $j$-the place.
Anti-symmetrizing the $\alpha$ and $\beta$ indexes, we
 obtain
\begin{equation}\label{eq1}
\Delta*\Omega +*\Omega |A|^2+2\sum_{\alpha<\beta,
i<j }\Omega_{\alpha \beta ij} R_{\alpha \beta i j}=0\\
\end{equation}
 This
formula essentially appeared in \cite{fc}. The parabolic analogue
was rediscovered by the author in the study of mean curvature
flows in higher codimension.

Another basic equation is equation (7.2) in \cite{wa1}.

\begin{equation}\label{eq2}\Delta |A|^2=2|\nabla A|^2- 2\sum_{i,j, m,k}(\sum_\alpha
h_{\alpha ij} h_{\alpha mk})^2-2\sum_{\alpha, \beta, i,j}(
R_{\alpha \beta ij})^2
\end{equation}

In codimension one case, this is the so called Simon's identity
which has been enormously useful in the study of minimal
hypersurfaces.

Next, we generalize a Lemma of \cite{ssy} to higher codimension.


\begin{lem}
\[|\nabla A|^2-|\nabla|A||^2\geq \frac{2}{n}|\nabla|A||^2\]
\end{lem}

\begin{proof}
\[
|\nabla A|^2-|\nabla|A||^2=\sum_{\alpha i j k }h_{\alpha ij,
k}^2-|A|^{-2}\sum_{k}(\sum_{\alpha i j} h_{\alpha ij}h_{\alpha ij,
k})^2\]

By expanding the right hand side, it is not hard to see

\[\begin{split}|\nabla A|^2-|\nabla |A||^2&=\frac{1}{2 |A|^2}\sum_{\alpha \beta
ijk rs}(h_{\alpha ij}h_{\beta rs, k}-h_{\beta rs} h_{\alpha
ij,k})^2\\
\end{split}\]

Recall the $h_{\alpha ij}$ are simultaneously diagonalizable.

So the expression is equal to \[\begin{split} &\sum_{\alpha i j k
}h_{\alpha ij, k}^2-|A|^{-2}\sum_{k}(\sum_{\alpha i} h_{\alpha
ii}h_{\alpha ii, k})^2\\
&\geq  \sum_{\alpha i j k }h_{\alpha ij, k}^2-|A|^{-2}\sum_{\alpha
i} h_{\alpha
ii}^2 \sum_{\alpha i k}h_{\alpha ii, k}^2\\
&\geq\sum_{\alpha, i\not= j ,k} h_{\alpha ij, k}^2 \\
&\geq \sum_{\alpha, i\not= j} h_{\alpha ij, i}^2+\sum_{\alpha,
i\not= j}
h_{ij, j}^2\\
&=2\sum_{\alpha, i\not= j}h_{\alpha ii, j}^2.\\
\end{split}\]

On the other hand, by diagonalization, we have

\[|\nabla |A||^2\leq \sum_{\alpha, ik} h_{\alpha ii, k}^2=\sum_{\alpha, i\not= k}
h_{\alpha ii,k}^2+\sum_i h_{\alpha ii,i}^2\]

By the minimal surface equation $\sum_i h_{\alpha ii}=0$, thus

\[|\nabla |A||^2\leq \sum_{\alpha, i\not= k} h_{\alpha ii,k}^2+\sum_i(\sum_{\alpha,
j\not= i} h_{\alpha jj,i})^2\leq \sum_{i\not=
j}h_{ii,j}^2+(n-1)\sum_{i\not= j} h_{\alpha ii,
j}^2=n\sum_{\alpha, i\not= j} h_{\alpha ii, j}^2
\]

\end{proof}

\section{Proofs of Theorem 1.1. and 1.2.}

{\em Proof of Theorem 1.1.} Since the normal bundle of $\Sigma$ is
flat, by (\ref{eq1}), we have

\[\Delta *\Omega+*\Omega |A|^2=0.\]

As $*\Omega>0$, this equation implies the first nonzero eigenvalue
of the operator $-\Delta-|A|^2$ is non-negative or that

\begin{equation}\label{stability}\int_\Sigma  u^2 |A|^2 \leq \int_\Sigma |\nabla u|^2\end{equation} for any $C^1$
function $u$. This follows from a well known argument, see for
example Lemma 1.24 (page 21) of \cite{cm}. Indeed, take the $\log$
of $*\Omega$ and compute

\[\Delta \log *\Omega=-|A|^2-|\nabla \log *\Omega|^2.\]

Multiply both sides by $u^2$, integrate by parts, apply the
Cauchy-Schwarz inequality, and the result is obtained.

We recall that for a minimal submanifold of $\R^{n+m}$, the
stability condition is equivalent to

\[\int_\Sigma \sum_{i, j}\langle \nabla_{e_i} e_j, V\rangle^2\leq
\int_\Sigma \sum_i |(\nabla_{e_i} V)^\perp|^2\] for any
compact-supported section $V$ of the normal bundle.

Since the normal bundle is flat, we can find parallel sections
$e_\alpha$ of the normal bundle and write $V=V^\alpha e_\alpha$.
Then the stability condition is equivalent to

\begin{equation}\label{sta}\int_\Sigma \sum_{i, j}(\sum_\alpha V^\alpha h_{\alpha ij})^2 \leq
\int_\Sigma \sum_\alpha |\nabla V^\alpha|^2\end{equation}

Apply (\ref{stability}) to $u=V^\alpha$ and sum over $\alpha$, we
derive
\begin{equation}\int_\Sigma \sum_\alpha (V^\alpha)^2 |A|^2 \leq
\int_\Sigma \sum_\alpha |\nabla V^\alpha|^2\end{equation}

This clearly implies (\ref{sta}).

\hfill $\Box$

Before proving Theorem 1.2, we first generalize an integral
curvature estimate of \cite{ssy}.
\begin{thm}
Let $\Sigma^n$ be a minimal graph with flat normal bundle, for
$p\in [2, 2+\sqrt{2/n})$ and $\phi>0$, $\phi \in C^1_c(\Sigma)$,
we have

\[\int |A|^{2p}\phi^{2p} \leq C(n, p)\int |\nabla \phi|^{2p}.\]
\end{thm}

\begin{proof} Set $u=|A|^{1+q} f$ in the stability inequality (\ref{stability}),
we obtain

\begin{equation}\label{stab}\begin{split}&\int |A|^{4+2q} f^2\\
&\leq \int |f\nabla |A|^{1+q}+|A|^{1+q} \nabla f|^2\\
&=(1+q)^2\int f^2|\nabla |A||^2|A|^{2q}+\int |A|^{2+2q}|\nabla
f|^2+ 2(1+q)\int f|A|^{1+2q}\nabla f\cdot \nabla
|A|\end{split}\end{equation}

We shall estimate the first term using the inequality

\[|A|\Delta |A|+|A|^4\geq \frac{2}{n}|\nabla |A||^2\]
which follow from Lemma 2.1 and equation (2.3).

Multiply both sides by $|A|^{2q} f^2$ and integrate by parts, we
derive

\begin{equation}\label{simon}\begin{split}&\frac{2}{n} \int |\nabla |A||^2 |A|^{2q} f^2\\&\leq \int
|A|^{4+2q}f^2-2\int f|A|^{1+2q} \nabla f\cdot \nabla
|A|^2-(1+2q)\int f^2|A|^{2q} |\nabla
|A||^2\end{split}\end{equation}

Substitute the inequality (\ref{stab}) for $\int |A|^{4+2q} f^2$
into (\ref{simon}), we obtain

\[(\frac{2}{n}-q^2)\int |A|^{2q} |\nabla |A||^2 f^2\leq \int
|A|^{2+2q} |\nabla f|^2+2q\int f|A|^{1+2q}\nabla f\cdot \nabla
|A|^2\]

Using the inequality $2xy\leq \epsilon x^2+\frac{1}{\epsilon}y^2$,
we arrive at

\begin{equation}\label{*}(\frac{2}{n}-q^2-\epsilon q)\int f^2|A|^{2q} |\nabla
|A||^2\leq (1+\frac{q}{\epsilon})\int |\nabla
f|^2|A|^{2+2q}\end{equation}

We go back to (\ref{stab}) and apply $(x+y)^2\leq 2x^2+2y^2$ and
obtain
\begin{equation}\begin{split}&\int |A|^{4+2q} f^2\leq 2(1+q)^2\int f^2|\nabla |A||^2|A|^{2q}+2\int
|A|^{2+2q}|\nabla f|^2.\end{split}\end{equation}

In view of (\ref{*}), if $\frac{2}{n}-q^2-\epsilon q>0$, we have
\[\int |A|^{4+2q} f^2\leq C \int |A|^{2+2q}|\nabla f|^2\] for some
constant $C$.

Take $p=2+q$ and $f=\phi^p$ and use $xy\leq
\frac{x^a}{a}+\frac{y^b}{b}$ for $\frac{1}{a}+\frac{1}{b}=1$, we
achieve

\[\int |A|^{2p} \phi^{2p}\leq C(n, p)\int |\nabla \phi|^{2p}\]

The condition $\frac{2}{n}-q^2>0$ translates to $2\leq p<
2+\sqrt{2/n}$ in view of $p=2+q$.

\end{proof}

Take $\phi$ to be the standard cut-off function supported in $B_R$
and $\equiv 1$ in $B_{R/2}$ , we have

\begin{equation}\label{A^2p}\int_{\Sigma\cap B_{R/2}} |A|^{2p}\leq C(n,
p)R^{-2p}Vol(\Sigma \cap B_R)\end{equation} for $p\in [2,
2+\sqrt{2/n})$.

We recall the following sub-mean-value inequality for minimal
submanifolds (\cite{bdm}, \cite{ms}).

\begin{thm} If $\Delta u\geq -Q u$ in $\Sigma\cap B_R$, $u\geq 0$, $Q\geq 0$ and
$Q\in L^q$ for some $q>n/2$ then

\[\sup_{\Sigma\cap B_{R/2}} u\leq C(R^{-n}\int_{\Sigma \cap B_R} u^2)^{1/2}\]where $C$ is
 a constant depending on $n$, $p$, $R^{2q-n}\int_{B_R} Q^q$ and the
 isoparametric constant.
\end{thm}

\noindent{\em Proof of Theorem 1.2.} We recall from (\ref{eq2})
$|A|^2$ satisfies $\Delta |A|^2+2|A|^4\geq 0$ take $u=|A|^2$ and
$Q=2|A|^2$, we have

\[\sup_{B_{R/2}} |A|^2\leq C(R^{-n}\int_{B_R} |A|^4)^{1/2}\]for
some $C$ that depends on $R^{2q-n}\int_{B_R} |A|^{2q}$. This
quantity is bounded by (3.8) and finite density assumption.

Suppose there exists a $q$ satisfying $2\leq q< 2+\sqrt{2/n}$ and
$q>\frac{n}{2}$. Take $p-2$ in Theorem 3.1, we have
\[\int_{\Sigma\cap B_{R/2}}|A|^4\leq C R^{-4}Vol(\Sigma\cap B_R).\]

We have \[\sup_{B_{R/2}} |A|^2\leq C(R^{-(n+4)}Vol(\Sigma\cap
B_R))^{1/2}\] Applying the assumption $Vol(\Sigma \cap B_R)\leq
c(n) R^{n}$ and let $R\rightarrow \infty$, we obtain $|A|^2=0$.
The relation satisfied by $q$ requires that

\[n<4+\sqrt{8/n}\] or $n=1\cdots 5$

\section{Proof of Theorem 1.3}

\noindent {\em Proof of Theorem 1.3.} We follow the proof of
Ecker-Huisken \cite{eh} closely. Since $R_{\alpha\beta ij}=0$, by
(\ref{eq1}) and (\ref{eq2}), we obtain

\begin{equation}\label{eq11}
\Delta*\Omega + |A|^2 *\Omega=0\\
\end{equation}
and
\begin{equation}\label{eq21}\Delta |A|^2=2|\nabla A|^2- 2\sum_{i,j, m,k}(\sum_\alpha
h_{\alpha ij} h_{\alpha mk})^2\geq 2|\nabla A|^2- 2|A|^4
\end{equation}

These two equations correspond to equation(2) and equation (3) in
\cite{eh}. $R_{\alpha\beta ij}=0 $ implies the matrices
$A^\alpha=\begin{bmatrix} h_{\alpha ij}\end{bmatrix}$ are pairwise
commutative and thus simultaneously diagonalizable. As in
\cite{ssy}, we can show
\[|\nabla A|^2\geq (1+\frac{2}{n})|\nabla |A||^2.\]
This is then identical to equation (4) in \cite{eh}.

As in \cite{eh}, for $p\geq 2$ and $q(1-2/n)\leq p-1+2/n$, we
derive

\[\Delta(|A|^p *\Omega^{-q})\geq(q-p) |A|^{p+2} *\Omega^{-q}.\]
Choose $q=p\geq (n-1)/2$, we obtain

\[\Delta (|A|^pv^p)\geq 0.\]

The sub-mean value inequality for subharmonic functions together
with the volume growth assumption yield

\begin{equation}\label{submin}|A|^pv^p(0)\leq c(n) R^{-n/2}\left(\int_{\Sigma\cap B_R}
|A|^{2p} v^{2p}\right)^{1/2}.\end{equation}

On the other hand, for $p\geq \max(3, n-1)$ fixed,

\[\Delta (|A|^{p-1} v^p)\geq |A|^{p+1}v^p\]

Multiply this equation by $|A|^{p-1} v^p \eta^{2p}$ where $\eta$
is  a test function and integrate by parts, we arrive at

\[\int_\Sigma |A|^{2p} v^{2p} \eta^{2p}\leq c(p)\int_\Sigma v^{2p}
|\nabla \eta|^{2p}.\] Take $\eta$ to be the standard cut-off
function such that $\eta\equiv 1$ on $B_R$, $\eta\equiv 0$ outside
$B_{2R}$ and $|\nabla \eta|\leq \frac{2}{R}$. Combine this with
\ref{submin}, apply the growth on the volume and $v$, and let
$R\rightarrow \infty$, we see that $|A|\equiv 0$. \hfill $\Box$

Comparing with Ecker-Huisken's proof, we notice the arguments only
differ in that we need to make the assumption of the growth of the
volume. A minimal graph of codimension one is area-minimizing and
a comparison argument gives this area bound. However, in the
higher codimensional case, we only prove the stability. It is
interesting to investigate whether a minimal graph with flat
normal bundle is area-minimizing.

\end{document}